\documentclass[12pt]{article}
\usepackage{amsmath,amsthm,amssymb}
\usepackage{times}
\usepackage{enumerate}
\usepackage{amsfonts}
\usepackage{fancyhdr}
\usepackage{titlesec}
\usepackage{cite}
\usepackage{ifthen}
\usepackage{amssymb}
\usepackage{fancyhdr}
\usepackage{titlesec}
\usepackage[arrow,matrix]{xy}
\usepackage{pifont}
\usepackage{stmaryrd}
\usepackage{setspace}
\usepackage{indentfirst}
\usepackage{amsmath,amssymb,amscd,bbm,amsthm,mathrsfs,dsfont}
\theoremstyle{definition}
\newtheorem{theorem}{Theorem}[section]

\newtheorem{claim}[theorem]{Claim}

\numberwithin{equation}{section}

\newcommand{\subjclass}[1]{\bigskip\noindent\emph{2010 Mathematics Subject Classification:}\enspace#1}
\newcommand{\keywords}[1]{\noindent\emph{Keywords:}\enspace#1}

\begin{document}


\baselineskip=16pt


\title{Weak solutions  of the convective Cahn-Hilliard equation with degenerate mobility\thanks{Correspondence to: Dr. Xiaopeng Zhao, E-mail address: zhaoxiaopeng@jiangnan.edu.cn}}

\author{  Xiaopeng Zhao    \\
\small{ School of Science, Jiangnan University, Wuxi 214122, China}}

\date{}

\maketitle


\begin{abstract}
In this paper, the existence of weak solutions of a convective Cahn-Hilliard equation with degenerate mobility is studied. We first define a notion of weak solutions and establish a regularized problems. The existence of such solutions is obtained by considered the limits of the regularized problems.


\end{abstract}
\section{Introduction} \label{sect1}

The Cahn-Hilliard equation
\begin{equation} \label{1-1}
\left\{ \begin{aligned}
       \partial_tu=&\nabla\cdot \left(M\nabla\mu\right), \\
                  \mu:=&-\gamma\Delta u+\varphi(u),
                          \end{aligned} \right.
                          \end{equation}
which arises in the study of phase separation on cooling binary solutions such as glasses, alloys and polymer mixtures(see \cite{Cahn,NCS,NC}) is a widely used phenomenological diffuse-interface model. Here, $u(x,t)$ is the relative concentration of the two phases, $\Phi(u)=\int_0^u\varphi(s)ds$ is a double-well potential with two equal minima at $u^-<u^+$ corresponding to the two pure phases, $M$ is the diffusion mobility.
During the past years, there are many classical papers related to the study of Cahn-Hilliard equation, see e.g. Elliott and Zheng\cite{Zheng}, Dlotko\cite{Dlotko}, Temam\cite{Temam},   Gilardi, Miranville and Schimperna\cite{Gi}, Yin\cite{Y1}, Dai and Du\cite{D1} and the reference therein.

Since Cahn-Hilliard equation is only a phenomenological model, various modifications of it has been proposed in order to capture the dynamical picture of the phase transition phenomena better. To name only a few, viscous Cahn-Hilliard equation (see\cite{VCH2,vch}), Cahn-Hilliard-Gurtin equation (see\cite{CHG2,chg}),  Cahn-Hilliard-Cook equation (see\cite{CHC1}), Cahn-Hilliard-Brinkman equation (see\cite{CHB2,CHB3}) and so on.

Recently, more and more people are interested in the convective Cahn-Hilliard equation
\begin{equation} \label{1-2}
\left\{ \begin{aligned}
       \partial_tu=&\nabla\cdot \left(M\nabla\mu\right)+\beta\cdot\nabla\psi(u), \\
                  \mu:=&-\gamma\Delta u+\varphi(u),
                          \end{aligned} \right.
                          \end{equation}
 which arises
naturally as a continuous model for the formation of facets and
corners in crystal growth (see \cite{GAA,WSJ}). Here, $u(x,t)$ denotes the
slope of the interface. The convective term
$\beta\cdot\nabla \psi(u)$ ( see \cite{GAA}) stems
from the effect of kinetic that provides an independent flux of the
order parameter, similar to the effect of an external field in
spinodal decomposition of a driven system.
In\cite{ZMPNG}, Zaks et al.
investigated the bifurcations of stations periodic solutions of a
convective Cahn-Hilliard equation; Eden and Kalantarov \cite{Eden1,Eden2} established
some results on the existence of a compact attractor for the
convective Cahn-Hilliard equation with periodic boundary conditions
in one space dimension and three space dimension; Zhao et al.\cite{Zhao2} considered the global attractor for the convective Cahn-Hilliard equation in two space dimension by using an iteration procedure and regularity estimates for the linear
semigroups. In \cite{Zhao}, the existence of optimal solutions and optimality condition for 2D convective Cahn-Hilliard equation were studied by Zhao and Liu.

As we all know, the mobility of many types of higher order nonlinear diffusion equations  are non-constant.
Recently, more and more authors paid their attentions to the well-posedness of solutions for higher-order diffusion equations with non-constant mobility (see e.g. Elliott and Garcke\cite{Elliott}, Schimperna\cite{Sch},
Bonetti, Dreyer and Schimperna\cite{Bonetti},  Schimperna and Zelik\cite{SZ} and
the references therein).

Since the main tool to study the higher order nonlinear diffusion equation with non-constant mobility is energy estimate. But for the convective Cahn-Hilliard equation together with non-constant mobility, we can't find the  free energy functional. This increase the difficulty of study. So far by now, there are only a few papers related to  the convective Cahn-Hilliard equation with non-constant mobility. In\cite{LY}, by using Schauder type estimates, Liu and Yin study the global existence of classical solutions for such equation with concentration dependent mobility. Moreover,
for the convective Cahn-Hilliard with degenerate mobility in 1D case, Liu\cite{LCC}
studied the global existence,  nonnegativity and the finite speed of
propagation of perturbations of solutions.

In this paper, we investigate the  convective Cahn-Hilliard equation
\begin{equation} \label{1-3}
\left\{ \begin{aligned}
       \partial_tu=&\nabla\cdot \left(M(u)\nabla\mu\right)+\beta\cdot\nabla\psi(u),\quad\forall x\in\Omega\subset\mathbb{R}^n,~~t\in[0,\infty),  \\
                  \mu:=&-\gamma\Delta u+\varphi(u),
                          \end{aligned} \right.
                          \end{equation}
 together with periodic boundary value conditions and initial value condition
 \begin{equation} \label{1-4}
u(x,0)=u_0(x),\quad x\in\Omega.
 \end{equation}
 For convenience, we suppose the mobility
 \begin{equation}
 \label{1-5}
 M(u)=|u|^{2m},\quad\hbox{for~all}~u\in\mathbb{R},
 \end{equation}where $m$ can be any positive number $0<m<\infty$ if the dimension $n=1,2$ and we require $0<m<\frac{n-1}{n-2}$ if $n\geq3$.
Moreover, suppose the polynomial forms of $\varphi(\chi)$ and $\psi(\chi)$  as the particular cases:
\begin{equation}
\label{1-6}
\varphi(\chi)=\sum_{i=1}^{2k+1}a_i\chi^i\quad\hbox{with}~a_i\in\mathbb{R},~a_{2k+1}>0, ~k\in \mathbb{Z}^+.
\end{equation}
and\begin{equation}
\label{1-7}
\psi(\chi)=\sum_{i=0}^{k}b_i\chi^{i+m}\quad\hbox{with}~b_i\in\mathbb{R},
\end{equation}
where $0\leq k<\infty$ if $n=1,2$ and $0\leq k\leq\frac1{n-2}$ if $n\geq3$.

For equation (\ref{1-3})-(\ref{1-4}), since the mobility is a degenerate function, the problem does not admit classical solutions in general. Our analysis involves two steps.

First, we approximate the degenerate mobility $M(u)=|u|^{2m}$ $(0<m<1)$ by a non-degenerate function $M_{\theta}(u)$ defined for a $\theta>0$ by
\begin{equation} \label{1-8}
M_{\theta}(u)=\left\{ \begin{aligned}
       |u|^{2m},\quad&\hbox{if}~|u|^2>\theta,\\
       \theta^m,\quad&\hbox{if}~|u|^2\leq\theta.
                          \end{aligned} \right.
                          \end{equation}
Since there exists a  uniform lower bound of $M_{\theta}(u)$, then we can find a sufficiently regular weak solution for (\ref{1-3})-(\ref{1-4}) with a mobility $M_{\theta}(u)$.

\begin{theorem}
\label{thm1.1}For any $u_0\in H^1(\Omega)$ and $T>0$ there is a function $u_{\theta}(x,t)$ such that
\begin{itemize}
\item $u_{\theta}\in L^{\infty}(0,T;H^1(\Omega))\bigcap C([0,T];L^p(\Omega))\bigcap L^2(0,T;H^3(\Omega))$, where $1\leq p<\infty$ if $n=1,2$ and $1\leq p<\frac{2n}{n-2}$ if $n\geq3$;
\item $\partial_t u_{\theta}\in L^2(0,T;(H^2(\Omega)')$;
\item $u_{\theta}(x,0)=u_0(x)$ for all $x\in\Omega$ which satisfies the following weak sense
\begin{equation}
\begin{aligned}
&\int_0^T\langle\partial_tu_{\theta},\rho\rangle_{((H^2(\Omega))',H^2(\Omega))}dt
\\
=&-\int_0^T\int_{\Omega}M_{\theta}(u_{\theta})\nabla\left[-\gamma\Delta u_{\theta}+\varphi(u_{\theta})\right]\cdot\nabla\rho dxdt-\beta\cdot\int_0^T\int_{\Omega}\psi(u_{\theta})\nabla\rho dxdt,
\label{1-9}
\end{aligned}\end{equation}for all $\rho\in L^2(0,T;H^2(\Omega))$.

\end{itemize}
\end{theorem}

Second, we consider the limit of $u_{\theta}$ as $\theta\rightarrow0$.  The limiting value $u$, of the functions $u_{\theta}$, does exist and, in a weak sense, solves the convective Cahn-Hilliard equation (\ref{1-6}) coupled with the mobility $M(u)$. It means that $u$ solves the convective Cahn-Hilliard equation in any open set $U\subset\Omega_T:=\Omega\times(0,T)$ where $u$ has enough regularity, namely where $\nabla\Delta u\in L^q(U)$ for some $q>1$. As for the singular set where $\nabla\Delta u$ fails to satisfy such a regularity condition, the singular set is contained in the set where $M(u)$ is degenerate, plus another set of Lebesgue measure zero.

\begin{theorem}
\label{thm1.2}For any $u_0\in H^1(\Omega)$ and $T>0$ there is a function $u:\Omega_T\rightarrow\mathbb{R}$ satisfying
\begin{itemize}
\item $u\in L^{\infty}(0,T;H^1(\Omega))\bigcap C([0,T];L^p(\Omega))$, where $1\leq p<\infty$ if $n=1,2$ and $1\leq p<6$ if $n\geq3$;
\item $\partial_t u\in L^2(0,T;(H^2(\Omega)')$;
\item $u(x,0)=u_0(x)$ for all $x\in\Omega$ which  can be considered as a weak solution for the equation in the following weak sense:
\begin{itemize}
\item Define $P$ as the set where $M(u)$ is not degenerate, that is
$$
P:=\{(x,t)\in\Omega_T:|u|^2\neq0\}.$$
There exist a set $B\subset\Omega_T$ with $|\Omega_T\setminus B|=0$ and a function $\zeta:\Omega_T\rightarrow\mathbb{R}^n$ satisfying $\chi_{B\bigcap P}M(u)\zeta\in L^2(0,T;L^{\frac{2n}{2+n}}(\Omega;\mathbb{R}^n))$, where $\chi_{B\bigcap P}$ is the characteristic function of $B\bigcap P$, such that
\begin{equation}\begin{aligned}
\int_0^T\langle\partial_tu,\rho\rangle_{((H^2(\Omega))',H^2(\Omega))}dt=-\int_{B\bigcap P}M(u)\zeta\cdot\nabla\rho dxdt-\beta\cdot\int_{\Omega}\psi(u)\nabla\rho dxdt,
\end{aligned}\nonumber\end{equation}for all $\psi\in L^2(0,T;H^2(\Omega))$;
\item Suppose that $\nabla\Delta u$ is the generalized derivative of $u$ in the sense of distributions. If $\nabla\Delta u\in L^q(U)$ for some open subset $U\subset\Omega_T$ and some $q>1$, then
    \begin{equation}\begin{aligned}
\zeta=&-\gamma\nabla\Delta u+\varphi'(u)\nabla u\quad\hbox{in}~U.
\end{aligned}\label{1-10}\end{equation}

\end{itemize}
\end{itemize}
\end{theorem}



The purpose of this paper is to study the existence of  weak solutions for problem (\ref{1-6})-(\ref{1-8}) with degenerate mobility $M(u)=|u|^{2m}$, where  $m$ can be any positive number $0<m<\infty$ if the dimension $n=1,2$ and  $0<m<\frac{n-1}{n-2}$ if $n\geq3$.
The main features and difficulties of this equation are caused by the lack of free energy functional,    the lack of maximum principle, a strong constraint imposed by the presence of the nonlinear principal part and the degeneracy.
Due to the lack of free energy, we have to resort to different techniques and modify the formulation accordingly. Our method is based on uniform
a prior estimates for local in time solutions. Because of the degeneracy, we  first consider the existence of weak solutions for the regularized problem, then based on the uniform estimates for the approximate solutions, we obtain the existence of weak solutions for the equation.

The outline of this paper is as follows. 
In Section 2, the existence of weak solutions for convective Cahn-Hilliard equation with positive mobilities is  studied. In Section 3, the existence of weak solutions for the equation with degenerate mobility is proved by considering the limits of such equation with positive mobilities. Finally in Section 4, conclusions are obtained.

Throughout this paper, the same letters $C$ and $C_i$ $(i=0,1,2,\cdots)$ denote positive constants that may depend on $T$, $\Omega$, $u_0$ but nothing else.

\section{Equation with positive mobility}

In this section, we consider the existence of weak solutions for convective Cahn-Hilliard equation with positive mobilities. For simplicity, we take $\Omega=[0,2\pi]^n$.  Write $\mathbb{Z}_+$ as the set of nonnegative integers. Then
$$
S_N=\{(2\pi)^{-\frac n2},\hbox{Re}(\pi^{-\frac n2}e^{i\xi\cdot x}),\hbox{Im}(\pi^{-\frac n2}e^{i\xi\cdot x}):\xi\in\mathbb{Z}_+^n\setminus \{(0,0,\cdots,0)\}\},
$$
form a complete orthonormal basic for $L^2(\Omega)$ that are also orthogonal in $H^k(\Omega)$ for any $k\geq1$. We label the basis as $\{\rho_j:j=1,2,\cdots\}$ with $\rho_1=(2\pi)^{-\frac n2}$.

\subsection{Galerkin approximation}
In order to study the weak solution of the equation with positive mobility, we define $$u^N(x,t)=\sum_{j=1}^Nc_j^N(t)\rho_j(x),\quad\mu^N(x,t)=\sum_{j=1}^Nd_j^N(t)\rho_j(x).
$$
Then,  the following system of equations for $j=1,2,\cdots,N$ is hold:
\begin{equation}
\label{2-1}
\int_{\Omega}\partial u^N\rho_jdx=-\int_{\Omega}M_{\theta}(u^N)\nabla\mu^N\cdot\nabla\rho_j dx+\beta\cdot\int_{\Omega}\nabla\psi_j\rho_jdx,
\end{equation}
\begin{equation}
\label{2-2}
\int_{\Omega}\mu^N\rho_jdx=\int_{\Omega}\left[\gamma\nabla  u^N\nabla\rho_j+\varphi(u^N)\rho_j
\right]dx,
\end{equation}\begin{equation}
\label{2-3}
u^N(x,0)=\sum_{j=1}^N\left(\int_{\Omega}u_0\rho_jdx\right)\rho_j(x).
\end{equation}
In fact, (\ref{2-1})-(\ref{2-2}) can be rewritten, equivalently, as
\begin{equation}
\label{2-4}\int_{\Omega}\partial u^N\rho_jdx=\int_{\Omega}M_{\theta}(u^N)(\gamma \nabla\Delta u^N-\nabla\varphi(u^N))\nabla\rho_jdx-\beta\cdot\int_{\Omega}\psi_j\nabla\rho_jdx.
\end{equation}

It is easy to check that problem (\ref{2-1})-(\ref{2-3}) (or problem (\ref{2-3})-(\ref{2-4})) is a system of ordinary differential equations for $\{c_j^N(t)\}_{j=1}^N$. Note that the right hand side of (\ref{2-1}) is continuous in $c_j$, the system has a local solution.

Setting $\rho_j=\gamma\Delta u^N-\varphi(u^N)$, we derive that
\begin{equation}
\nonumber\begin{aligned}&
\int_{\Omega}\partial_tu^N(\gamma\Delta u^N-\varphi(u^N))dx+\int_{\Omega}\nabla\cdot[M_{\theta}(u_{\theta})\nabla(\gamma\Delta u^N-\varphi(u^N))](\gamma\Delta u^N-\varphi(u^N))dx
\\
=&-\beta\cdot\int_{\Omega}\psi(u^N)\nabla(\gamma\Delta u^N-\varphi(u^N))dx.
\end{aligned}\end{equation}
By Young's inequality, we deduce that
\begin{equation}\nonumber\begin{aligned}
&\frac d{dt}\int_{\Omega}\left(\frac{\gamma}2|\nabla u^N|^2+\Phi(u^N)\right)dx+\int_{\Omega}M_{\theta}(u_{\theta})|\gamma\nabla\Delta u^N-\nabla\varphi(u^N)|^2dx
\\
\leq&\frac12\int_{\Omega}M_{\theta}(u_{\theta})|\gamma\nabla\Delta u^N-\nabla\varphi(u^N)|^2dx+\frac{|\beta|^2}2\int_{\Omega}\frac{|\psi(u^N)|^2}{M_{\theta}(u^N)}dx,
\end{aligned}\end{equation}
that is
\begin{equation}
\begin{aligned}&\frac d{dt}\int_{\Omega}\left(\frac{\gamma}2|\nabla u^N|^2+\Phi(u^N)\right)dx+\frac12\int_{\Omega}M_{\theta}(u_{\theta})|\gamma\nabla\Delta u^N-\nabla\varphi(u^N)|^2dx
\\
\leq&\frac{|\beta|^2}2\left(C_1\int_{\Omega}\Phi(u^N)dx+C_2\right)\leq\frac{|\beta|^2}2\left[C_1\int_{\Omega}\left(\frac{\gamma}2|\nabla u^N|^2+\Phi(u^N)\right)dx+C_2\right].
\label{2-5}
\end{aligned}\end{equation}
Using Gronwall's inequlity, we obtain
\begin{equation}\begin{aligned}
\label{2-6}&
\int_{\Omega}\left(\frac{\gamma}2|\nabla u^N|^2+\Phi(u^N)\right)dx\\\leq
&e^{\frac{C_1|\beta|^2t}2}\left[\int_{\Omega}\left(\frac{\gamma}2|\nabla u^N(x,0)|^2+\Phi(u^N(x,0))\right)dx+C_3\right]
\\
\leq&\frac{\gamma}2e^{\frac{C_1|\beta|^2T}2}\|\nabla u_0\|^2+e^{\frac{C_1|\beta|^2T}2}\int_{\Omega}\Phi(u_0)dx+C_3e^{\frac{C_1|\beta|^2T}2}
\\
\leq&C(\|\nabla u_0\|^2+\int_{\Omega}|u_0|^{2k+1}dx+1)
\leq C(\|\nabla u_0\|^2+\|u_0\|_{H^1}^{2k+1}+1)\leq C.
\end{aligned}\end{equation}

Taking $j=1$ in (\ref{2-1}), we have $\int_{\Omega}\partial_tu^Ndx=0$. Then,
$$
\int_{\Omega}u^N(x,t)dx=\int_{\Omega}u^N(x,0)dx.
$$
We define $\Pi_N$ as the $L^2$ projection operator from $L^2(\Omega)$ into $\hbox{span}\{\rho_j\}_{j=1}^N$, which means $\Pi_N\rho:=\sum_{j=1}^N\left(\int_{\Omega}\rho\rho_jdx\right)\rho_j.
$
 Hence,
\begin{equation}
\label{2-7}
\left|\int_{\Omega}u^N(x,t)dx\right|=\left|\int_{\Omega}u^N(x,0)dx\right|\leq\|u_0\||\Omega|^{\frac12}.
\end{equation}
It then follows from (\ref{2-6}) and (\ref{2-7}) that
\begin{equation}
\label{2-8}
\|u^N\|_{L^2(0,T;H^1(\Omega))}\leq C,\quad \hbox{for~all}~N.
\end{equation}
Integrating (\ref{2-5}) with respect to $t$ over $(0,T)$, we have
\begin{equation}
\label{2-9}
\|\sqrt{M_{\theta}(u^N)}\nabla\mu^N\|_{L^2(\Omega_T)}\leq C,\quad \hbox{for~all}~N.
\end{equation}
Using Sobolev's embedding theorem, we get
\begin{equation}
\label{2-10}
\|\varphi(u^N)\|_{L^{\infty}(0,T;L^2(\Omega))}+\|\psi(u^N)\|_{L^{\infty}(0,T;L^2(\Omega))}\leq C,
\end{equation}
and
\begin{equation}
\label{2-11}
\|M_{\theta}(u^N)\|_{L^{\infty}(0,T;L^{\frac n2}(\Omega))}\leq C.
\end{equation}
Follows form (\ref{2-18}), the coefficients $\{c_j^N:j=1,2,\cdots,N\}$ are bounded in time and there is a global solution for the system (\ref{2-1})-(\ref{2-3}).

\subsection{Convergence of $\{u^N\}$}
For all $\rho\in L^2(0,T;H^2(\Omega))$, we set $\Pi_N\rho(x,t)=\sum_{j=1}^Na_j(t)\rho_j(x)$. Then
\begin{equation}
\label{2-13}
\begin{aligned}&
\left|\int_{\Omega}\partial_t u^N\rho dx\right|=\left|\int_{\Omega}\partial_tu^N\Pi_N\rho dx\right|
\\
\leq&\left|\int_{\Omega}M_{\theta}(u^N)\nabla\mu^N\cdot\nabla\Pi_N\rho dx\right|+|\beta|\left|\int_{\Omega}\psi(u^N)\nabla\Pi_N\rho dx\right|
\\
\leq&\|\sqrt{M_{\theta}(u^N)}\|_{L^{n}(\Omega)}\|\sqrt{M_{\theta}(u^N)}\nabla\mu^N\|_{L^2(\Omega)}\|\nabla\Pi_N\rho\|_{L^{\frac{2n}{n-2}}(\Omega)}
+|\beta|\|\psi(u^N)\|_{L^{\frac{2n}{n+2}}(\Omega)}\|\nabla\Pi_N\rho\|_{L^{\frac{2n}{n-2}}(\Omega)}
\\
\leq &C\left(\|\sqrt{M_{\theta}(u^N)}\nabla\mu^N\|_{L^2(\Omega)}+\|\psi(u^N)\|_{L^{\frac{2n}{n+2}}(\Omega)}\right)\|\rho\|_{H^2(\Omega)}
\\
\leq&C\left(\|\sqrt{M_{\theta}(u^N)}\nabla\mu^N\|_{L^2(\Omega)}+\|u^N\|_{H^1(\Omega)}^{\frac{n+2}{2n-4}}+1\right)\|\rho\|_{H^2(\Omega)}.
\end{aligned}\end{equation}
By (\ref{2-11}) and (\ref{2-13}), we get
\begin{equation}
\label{2-14}
\begin{aligned}
\left|\int_0^T\!\!\!\int_{\Omega}\partial_tu^N\psi dxdt\right|
\leq & C\left(\int_0^T\!\!\!\int_{\Omega}M_{\theta}(u^N)|\nabla\mu^N|^2dxdt\right)^{\frac12}\left(\int_0^T\|\rho\|^2_{H^2(\Omega)}dt\right)^{\frac12}
\\&+C\left(\int_0^T(\|u^N\|_{H^1(\Omega)}^{\frac{n+2}{2n-4}}+1)dt\right)^{\frac12}\left(\int_0^T\|\rho\|^2_{H^2(\Omega)}dt\right)^{\frac12}
\\
\leq &C\|\rho\|_{L^2(0,T;H^2(\Omega))}.
\end{aligned}\end{equation}
Thus
\begin{equation}
\label{2-15}
\|\partial_tu^N\|_{L^2(0,T;(H^2(\Omega))')}\leq C,\quad\hbox{for~all}~N.
\end{equation}
Note that the embedding $H^1(\Omega)\hookrightarrow L^q(\Omega)$ is compact for $1\leq q<\infty$ if $n=1,2$ and $1\leq q<\frac{2n}{n-2}$ if $n\geq3$. By Aubin-Lions' Lemma (see in \cite{Simon1,Simon2}), the embeddings
$$
\{u\in L^2(0,T;H^1(\Omega)):\partial_tu\in L^2(0,T;(H^2(\Omega))')\}\hookrightarrow L^2(0,T;L^q(\Omega)),
$$
$$
\{u\in L^{\infty}(0,T;H^1(\Omega)):\partial_tu\in L^2(0,T;(H^2(\Omega))')\}\hookrightarrow C([0,T];L^q(\Omega)),$$
are compact for the values of $q$ indicated.
The above boundedness of $\{u^N\}$ and $\{\partial_tu^N\}$ enables us to find a subsequence, not relabeled, and $u_{\theta}\in L^{\infty}(0,T;H^1(\Omega))$ such that as $N\rightarrow\infty$,
\begin{equation}
\label{2-16}
u^N\rightharpoonup u_{\theta}\quad\hbox{weakly}*~\hbox{in}~L^{\infty}(0,T;H^1(\Omega)),
\end{equation}
\begin{equation}
\label{2-17}
u^N\rightharpoonup u_{\theta}\quad\hbox{strongly~in}~C([0,T];L^q(\Omega)),
\end{equation}
\begin{equation}
\label{2-18}
u^N\rightharpoonup u_{\theta}\quad\hbox{strongly~in}~L^2(0,T;L^q(\Omega))~\hbox{and~almost~everywhere~in}~\Omega_T,
\end{equation}
\begin{equation}
\label{2-19}
\partial_tu^N\rightharpoonup \partial_tu_{\theta}\quad\hbox{weakly~in}~L^2(0,T;(H^2(\Omega))'),
\end{equation}
where  $1\leq q<\infty$ if $n=1,2$ and $1\leq q<\frac{2n}{n-2}$ if $n\geq3$. Furthermore, we have
\begin{equation}\label{2-20}
\|u_{\theta}\|_{L^{\infty}(0,T;H^1(\Omega))}\leq C,
\end{equation}
\begin{equation}\label{2-21}
\|\partial_tu_{\theta}\|_{L^2(0,T;(H^2(\Omega))')}\leq C,
\end{equation}
Since $M_{\theta}(u^N)=|u^N|^{2m}$ is continuous and $\sqrt{M_{\theta}(u^N)}=|u^N|^m$. It then follows from (\ref{2-17}) and the general dominated convergence theorem (see, for example, Theorem 17 of Section 4.4 on p.92 of \cite{12}) that
\begin{equation}
\label{2-22}
M_{\theta}(u^N)\rightarrow M_{\theta}(u_{\theta}),\quad\hbox{strongly~in}~C([0,T];L^{\frac n2}(\Omega)),
\end{equation}
\begin{equation}\label{2-23}
\sqrt{M_{\theta}(u^N)}\rightarrow\sqrt{M_{\theta}(u_{\theta})} ,\quad\hbox{strongly~in}~C([0,T];L^{n}(\Omega)).
\end{equation}
Moreover, since $|\varphi(u^N)|\leq C(|u^N|^{2k+1}+1)$ for $0\leq k\leq\frac1{n-2}$ when $n\geq3$, we have
\begin{equation}\label{2-24}
\varphi(u^N)\rightarrow \varphi(u),\quad\hbox{strongly~in}~C([0,T];L^q(\Omega)),
\end{equation} where $1\leq q< \infty$ if $n=1,2$ and $1\leq q<\frac{2n}{n-2}$ if $n\geq3$.
Using (\ref{2-10}), we can find a $w\in L^{\infty}(0,T;L^2(\Omega))$ such that $\varphi(u^N)\rightharpoonup w$ weakly* in $L^{\infty}(0,T;L^2(\Omega))$. Combining with (\ref{2-24}), we get $w=\varphi(u_{\theta})$. Thus,
\begin{equation}
\label{2-26}
\varphi(u^N)\rightharpoonup \varphi(u_{\theta}),\quad\hbox{weakly} *~\hbox{in}~L^{\infty}(0,T;L^2(\Omega)).
\end{equation}We have $|\psi(u^N)|\leq C(1+|u^N|^{m+k})\leq C(1+|u^N|^{\frac n{n-2}})$ when $n\geq3$. Then, by the same method as above,  we obtain
\begin{equation}\label{2-27}
\psi(u^N)\rightarrow \psi(u),\quad\hbox{strongly~in}~C([0,T];L^q(\Omega)),
\end{equation} where $1\leq q< \infty$ if $n=1,2$ and $1\leq q<\frac{2n}{n-2}$ if $n\geq3$.
Using (\ref{2-10}), we can find a $\vartheta\in L^{\infty}(0,T;L^2(\Omega))$ such that $\psi(u^N)\rightharpoonup \vartheta$ weakly* in $L^{\infty}(0,T;L^2(\Omega))$. Combining with (\ref{2-27}), we get $\vartheta=\psi(u_{\theta})$. Thus,
\begin{equation}
\label{2-28}
\psi(u^N)\rightharpoonup \psi(u_{\theta}),\quad\hbox{weakly} *~\hbox{in}~L^{\infty}(0,T;L^2(\Omega)).
\end{equation}
\subsection{Weak solution}


We have $M_{\theta}(u^N)\geq\theta^m>0$. Then, follows from (\ref{2-5}), we derive that
$$
\|\nabla\mu^N\|_{L^2(\Omega_T)}\leq C\theta^{-\frac m2}.
$$
Setting $\psi_j=1$ in (\ref{2-2}), we deduce that
\begin{equation}
\left|\int_{\Omega}\mu^Ndx\right|\leq C<\infty.\nonumber
\end{equation}
By Poincar\'{e}'s inequality, we get
\begin{equation}
\label{2-28}\|\mu^N\|_{L^2(0,T;H^1(\Omega)}\leq C(\theta^{-\frac m2}+1).
\end{equation}
Then, we can find a subsequence of $\mu^N$, not relabeled, and
$$
\mu_{\theta}\in L^2(0,T;H^1(\Omega)),
$$
such that
\begin{equation}
\label{2-29}
\mu^N\rightharpoonup \mu_{\theta},\quad\hbox{weakly~in}~L^2(0,T;H^1(\Omega)).
\end{equation}
Combining (\ref{2-29}) and (\ref{2-23}) together gives
$$
\sqrt{M_{\theta}(u^N)}\nabla\mu^N\rightharpoonup \sqrt{M_{\theta}(u_{\theta})}\nabla\mu_{\theta},\quad\hbox{weakly~in}~L^2(0,T;L^{\frac{2n}{2+n}}(\Omega)).
$$
Follows from (\ref{2-11}), $\sqrt{M_{\theta}(u^N)}\nabla\mu^N$ is bounded in $L^2(\Omega_T)$. Then we can extract a further sequence, not relabeled, such that the above weak convergence can be improved
\begin{equation}\label{2d-29}
\sqrt{M_{\theta}(u^N)}\nabla\mu^N\rightharpoonup \sqrt{M_{\theta}(u_{\theta})}\nabla\mu_{\theta},\quad\hbox{weakly~in}~L^2(\Omega_T).
\end{equation}
Hence
\begin{equation}
\label{2-30}
\int_0^T\int_{\Omega}M_{\theta}(u_{\theta})|\nabla\mu_{\theta}|^2dxdt\leq C<\infty.
\end{equation}
 We have $M_{\theta}(u^N)\nabla\mu^N=\sqrt{M_{\theta}(u^N)}\sqrt{M_{\theta}(u^N)}\nabla \mu^N$, for any $\Psi\in L^2(0,T;L^{\frac{2n}{n-2}}(\Omega))$, we have
\begin{equation}
\begin{aligned}&
\left|\int_0^T\!\!\!\int_{\Omega}(M_{\theta}(u^N)\nabla \mu^N-M_{\theta}(u_{\theta})\nabla\mu_{\theta})\Psi dxdt\right|
\\
\leq&\int_0^T\|\sqrt{M_{\theta}(u^N)}-\sqrt{M_{\theta}(u_{\theta})}\|_{L^3(\Omega)}\|\sqrt{M_{\theta}(u^N)}\nabla\mu^N\|_{L^2(\Omega)}\|\Psi\|_{L^6(\Omega)}
\\
&+\left|\int_{\Omega_T}(\sqrt{M_{\theta}(u^N)}\nabla\mu^N-\sqrt{M_{\theta}(u_{\theta})}\nabla\mu_{\theta})\sqrt{M_{\theta}(u_{\theta})}\Psi dxdt\right|
\\
\leq&\sup_{t\in[0,T]}\|\sqrt{M_{\theta}(u^N)}-\sqrt{M_{\theta}(u_{\theta})}\|_{L^n(\Omega)}\|\sqrt{M_{\theta}(u^N)}\nabla\mu^N\|_{L^2(\Omega_T)}\|\Psi\|_{L^2(0,T;L^{\frac{2n}{n-2}}(\Omega))}
\\
&+\|\sqrt{M_{\theta}(u^N)}\nabla\mu^N-\sqrt{M_{\theta}(u_{\theta})}\nabla\mu_{\theta}\|_{L^2(\Omega_T)}\|M_{\theta}(u_{\theta})\|
_{L^{\infty}(0,T;L^{\frac n2}(\Omega))}\|\Psi\|_{L^2(0,T;L^{\frac{2n}{n-2}}(\Omega))}
\\
\rightarrow&0,\quad\hbox{as}~N\rightarrow\infty.
\end{aligned}\nonumber\end{equation}
Then
$$
M_{\theta}(u^N)\nabla\mu^N\rightharpoonup M_{\theta}(u_{\theta})\nabla\mu_{\theta},\quad\hbox{weakly~in}~L^2(0,T;L^{\frac{2n}{n-2}}(\Omega)).
$$
Then, for any $a(t)\in L^2(0,T)$, since $a(t)\nabla\psi_j\in L^2(0,T;L^{\frac{2n}{n-2}}(\Omega))$. Multiplying (\ref{2-1}) by $a(t)$, integrating in time over $(0,T)$, taking limit as $N\rightarrow\infty$, we obtain
\begin{equation}
\label{2-31}\begin{aligned}&
\int_0^T\langle\partial_tu_{\theta},a(t)\rho_j(x)\rangle_{((H^2(\Omega))',H^2(\Omega))}dt\\&=-\int_0^T\int_{\Omega}M_{\theta}(u_{\theta})\nabla\mu_{\theta}\cdot a(t)\nabla\rho_jdxdt-\beta\cdot\int_0^T\int_{\Omega}\psi_{\theta}a(t)\nabla\rho_j  dxdt,\quad\forall j\in\mathbb{N}.
\end{aligned}\end{equation}
For any $\rho\in L^2(0,T;H^2(\Omega))$, its Fourier series $\sum_{j=1}^{\infty}a_j(t)\rho_j$ converges strongly to $\rho\in L^2(0,T;H^2(\Omega))$. Then, $\sum_{j=1}^{\infty}a_j(t)\nabla\rho_j$ converges strongly to $\nabla\rho$ in $L^2(0,T;L^{\frac{2n}{n-2}}(\Omega))$,
\begin{equation}
\label{2-31}\begin{aligned}&
\int_0^T\langle\partial_tu_{\theta},\rho\rangle_{((H^2(\Omega))',H^2(\Omega))}dt\\&=-\int_0^T\!\!\!\int_{\Omega}M_{\theta}(u_{\theta})\nabla\mu_{\theta}\cdot \nabla\rho dxdt
-\beta\cdot\int_0^T\!\!\!\int_{\Omega}\psi_{\theta}\nabla\rho dxdt,~\forall \rho\in L^2(0,T;H^2(\Omega))
\end{aligned}\end{equation} Consider the initial value, we have
$$
u^N(x,0)\rightarrow u_{\theta}(x,0)\quad\hbox{as}~N\rightarrow\infty~\hbox{in}~L^2(\Omega).
$$
By (\ref{2-17}), $u_{\theta}(x,0)=u_0(x)$ in $L^2(\Omega)$.

\subsection{Regularity of $u_{\theta}$}
Note that $a_j(t)\rho_j\in L^2(0,T;C(\bar{\Omega}))$.
By (\ref{2-2}), for any $a_j(t)
\in L^2(0,T)$,  in the limit when $N\rightarrow\infty$, we get
\begin{equation}\begin{aligned}
\int_0^T\!\!\!\int_{\Omega}\mu_{\theta}a_j(t)\rho_jdxdt=\int_0^T\!\!\!\int_{\Omega}\left[\gamma\nabla u_{\theta}\cdot a_j(t)\nabla\rho_j+\varphi'(u_{\theta})a_j(t)\rho_j\right]dxdt,\end{aligned}
\end{equation}
for all $j\in\mathbb{N}$. Then, for any $\rho\in L^2(0,T;H^1(\Omega))$, since its Fourier series strongly converges to $\rho$ in $L^2(0,T;H^1(\Omega))$, we derive that
\begin{equation}\begin{aligned}
\nonumber \int_0^T\!\!\!\int_{\Omega}\mu_{\theta}\psi dxdt=\int_0^T\!\!\!\int_{\Omega}\left[\gamma\nabla u_{\theta}\cdot \nabla\rho+\varphi(u_{\theta})\rho\right]dxdt.\end{aligned}
\end{equation}
Note that $\varphi(u_{\theta})\in L^{\infty}(0,T;L^2(\Omega))$ and $\mu_{\theta}\in L^2(0,T;H^1(\Omega))$. By regularity theory, we see that $u_{\theta}\in L^2(0,T;H^2(\Omega))$. Hence
\begin{equation}\label{2-32}
\mu_{\theta}=\gamma\Delta u_{\theta}+\varphi(u_{\theta}),\quad\hbox{almost~everywhere~in}~\Omega_T.
\end{equation}
We have $|\varphi'(u_{\theta})|\leq C(1+|u_{\theta}|^{2k})$ and $u_{\theta}\in L^{\infty}(0,T;H^1(\Omega))\hookrightarrow L^{\infty}(0,T;L^q(\Omega))$. where $1\leq q<\infty$ if $n=1,2$ and $1\leq q<\frac{2n}{n-2}$ if $n\geq3$. Thus
\begin{equation}\begin{aligned}
\label{2-33}
\int_0^T\!\!\!\int_{\Omega}|\nabla \varphi(u_{\theta})|^2dxdt=&\int_0^T\!\!\!\int_{\Omega}|\varphi'(u_{\theta})|^2|\nabla u_{\theta}|^2dxdt
\leq C\left(\|u_{\theta}\|_{H^1(\Omega)}^{4k}+1\right)\|u_{\theta}\|^2_{L^2(0,T;H^2(\Omega))}.
\end{aligned}\end{equation}
Combining with $\mu_{\theta}\in L^2(0,T;H^1(\Omega))$, by (\ref{2-32}), we get $u_{\theta}\in L^2(0,T;H^3(\Omega))$ and
\begin{equation}
\label{2-37}
\nabla\mu_{\theta}=-\gamma\nabla\Delta u_{\theta}+\varphi'(u_{\theta})\nabla u_{\theta},\quad\hbox{almost~everywhere~in}~\Omega_T.
\end{equation}
%
It then follows from (\ref{2-31}) and (\ref{2-37}) that
\begin{equation}
\begin{aligned}
\label{2-38}&
\int_0^T\langle\partial_tu_{\theta},\psi\rangle_{((H^3(\Omega))',H^3(\Omega))}
\\=&-\int_0^T\!\!\!\int_{\Omega}M_{\theta}(u_{\theta})(-\gamma\nabla\Delta u_{\theta}+\varphi'(u_{\theta})\nabla u_{\theta})\cdot\nabla\rho dxdt-\beta\cdot\int_0^T\!\!\!\int_{\Omega}\psi(u_{\theta})\nabla \rho dxdt,
\end{aligned}\end{equation}for all $\psi\in L^2(0,T;H^2(\Omega))$.

\section{Equation with degenerate mobility}

In this section, we consider the convective Cahn-Hilliard equation with degenerate mobility.
The proof consists of two parts. First, the weak convergence of approximate solutions $u_{\theta_i}$ defined in the above section for a sequence of positive numbers $\theta_i\rightarrow0$ is proved; Second, the weak limit $u$ solves the degenerate equation in the weak sense is showed.

\subsection{Weak convergence of approximate solutions $\{u_{\theta_i}\}$}

Fix $u_0\in H^1(\Omega)$ and a sequence $\theta_i>0$ that monotonically decreases to $0$ as $i\rightarrow\infty$. Then, based on Theorem \ref{thm1.1}, for any $\theta_i>0$, we have a $u_i$ such that
\begin{equation}\begin{aligned}&
u_i\in L^{\infty}(0,T;H^1(\Omega))\bigcap C([0,T;L^q(\Omega))\bigcap L^2(0,T;H^3(\Omega)),\nonumber
\\&\partial_tu_i\in L^2(0,T;(H^2(\Omega))'),
\end{aligned}\end{equation}
where $1\leq q<\infty$ if $n=1,2$ and $1\leq p<\frac{2n}{n-2}$ if $n\geq3$, such that for all $\rho \in L^2(0,T;H^2(\Omega))$
\begin{equation}
\label{3-1}\begin{aligned}
\int_0^T\langle\partial_tu_i,\rho\rangle_{((H^2(\Omega))',H^2(\Omega))}dt=-\int_0^T\!\!\!\int_{\Omega}M_i(u_i)\nabla\mu_i\cdot\nabla\psi dxdt+\int_0^T\!\!\!\int_{\Omega}\nabla\psi(u_i)\rho dxdt,
\end{aligned}
\end{equation}
\begin{equation}
\label{3-2}
\mu_i=-\gamma\Delta u_i+\varphi(u_i),
\end{equation}
where $u_i=u_{\theta_i}$ and $M_i(u_i)=M_{\theta_i}(u_i)$. Based on the arguments in the above section, the bounds on the righthand side of (\ref{2-10}), 
 (\ref{2-15}), (\ref{2-30}) depend only on the growth conditions of the mobility. Hence, there exists a positive constant $C$ independent of $\theta_i$ such that
\begin{equation}
\label{3-3}
\|u_i\|_{L^{\infty}(0,T;H^1(\Omega))}\leq C,
\end{equation}
\begin{equation}
\label{3-4}
\|\partial_tu_i\|_{L^2(0,T;(H^2(\Omega))')}\leq C,
\end{equation}
\begin{equation}
\label{3-5}\|\sqrt{M_i(u_i)}\nabla\mu_i\|_{L^2(\Omega_T)}\leq C,
\end{equation}
\begin{equation}
\label{3-6}
\|\mu_i\|_{L^{\infty}(0,T;(H^1(\Omega))')}\leq C.
\end{equation}
\begin{equation}\label{3-7}
\|M_i(u_i)\|_{L^{\infty}(0,T;L^{\frac n2}(\Omega))}\leq C.
\end{equation}
\begin{equation}\label{3-7a}
\|\varphi(u_i)\|_{L^{\infty}(0,T;L^2(\Omega))}\leq C.
\end{equation}\begin{equation}\label{3-7b}
\|\psi(u_i)\|_{L^{\infty}(0,T;L^2(\Omega))}\leq C.
\end{equation}
Similar to the proof of Theorem 1.1, the above boundedness of $\{u_i\}$ and $\{\partial_tu_i\}$ enables us to find a subsequence, not relabeled, and $$u\in L^{\infty}(0,T;H^1(\Omega))\bigcap C([0,T];L^q(\Omega))$$ for any $1\leq q<\infty$ if $n=1,2$ and $1\leq q<\frac{2n}{n-2}$ if $n\geq3$, such that, as $i\rightarrow\infty$,
\begin{equation}
\label{3-8}
u_i\rightharpoonup u\quad \hbox{weakly} *~\hbox{in}~L^{\infty}(0,T;H^1(\Omega)),
\end{equation}
\begin{equation}\label{3-9}
u_i\rightarrow u\quad\hbox{strongly~in}~C([0,T];L^q(\Omega)),
\end{equation}
\begin{equation}
u_i\rightarrow u\quad\hbox{strongly~in}~L^2(0,T;L^q(\Omega))~\hbox{and~almost~everywhere~in}~\Omega_T,
\label{3-10}\end{equation}
\begin{equation}
\label{3-11}
\partial_tu_i\rightharpoonup \partial_tu\quad\hbox{weakly~in}~L^2(0,T;(H^2(\Omega))').\end{equation}
By general dominated convergence theorem, and the uniform convergence of $M_i\rightarrow M$ and $\sqrt{M_i}\rightarrow\sqrt{M}$ in $\mathbb{R}$ as $i\rightarrow\infty$, we have
\begin{equation}\label{3-12}
M_i(u_i)\rightarrow M(u)\quad\hbox{strongly~in}~C([0,T];L^{\frac n2}(\Omega)),
\end{equation}\begin{equation}\label{3-13}
\sqrt{M_i(u_i)}\rightarrow \sqrt{M(u)}\quad\hbox{strongly~in}~C([0,T];L^n(\Omega)).
\end{equation}
It then follows from (\ref{3-7b}) that
$$
\|\psi(u_i)\nabla u_i\|_{L^{\infty}(0,T;L^{\frac{2n}{2+n}}(\Omega))}\leq C.
$$
Then, we get
\begin{equation}
\label{nan}
\psi(u_i)\rightharpoonup \psi(u)\quad\hbox{weakly}-*\hbox{in}~L^{\infty}(0,T;L^{\frac{2n}{2+n}}(\Omega)).
\end{equation}
By (\ref{3-5}), there exists a $\xi\in L^2(\Omega_T)$ such that
\begin{equation}
\label{3-14}
\sqrt{M_i(u_i)}\nabla\mu_i\rightharpoonup\xi\quad\hbox{weakly~in}~L^2(\Omega_T).
\end{equation}
It then follows from (\ref{3-13})-(\ref{3-14})  that
\begin{equation}\label{3-15}
M_i(u_i)\nabla\mu_i\rightharpoonup\sqrt{M(u)}\xi\quad\hbox{weakly~in}~L^2(0,T;L^{\frac{2n}{2+n}}(\Omega)).
\end{equation}
Taking the limit as $i\rightarrow\infty$ in (\ref{3-1}), for all $\rho\in L^2(0,T;H^2(\Omega))$, we  get
\begin{equation}
\label{3-16}\begin{aligned}
\int_0^T\langle\partial_tu,\rho\rangle_{(H^2(\Omega),(H^2(\Omega))')}dt=-\int_0^T\!\!\!\int_{\Omega}\sqrt{M(u)}\xi\cdot\nabla\rho dxdt-\beta\cdot\int_0^T\!\!\!\int_{\Omega}\psi(u)\nabla \rho dxdt.\end{aligned}
\end{equation}
For the initial value, since $u_i(x,0)=u_0(x)$, using (\ref{3-9}), we derive that $u(x,0)=u_0(x)$.

\subsection{Weak solution for the degenerate equation}

We consider the relation between $\sqrt{M(u)}$ and $u(x,t)$. In fact, we need to consider the convergence properties of $\nabla\mu_i=-\gamma\nabla\Delta u_i+\varphi'(u_i)\nabla u_i$. Using the same method as \cite{D1}, we have
$$
\varphi'(u_i)\nabla u_i\rightharpoonup \varphi'(u_i)\nabla u\quad\hbox{weakly}*\hbox{in}~L^{\infty}(0,T;L^{\frac{2n}{2+n}}(\Omega)).
$$

Now, we consider the weak convergence of $\nabla\mu_i$ as a whole. Since the following discussion is so standard and can be performed essentially as in \cite{D1}, we will give very few details.

We choose a sequence of positive numbers
$\delta_j$ that monotonically decreases to $0$. By (\ref{3-10}) and Egorov's theorem, for every $\delta_j>0$, there exists a subset $B_j\subset\Omega_T$ with $|\Omega_T\backslash
B_j|<\delta_j$ such that
\begin{equation}
u_i\rightarrow u\quad\hbox{uniformly~in}B_j.\nonumber
\end{equation}
Take
\begin{equation}
\label{3-17}
B_1\subset B_2\subset\cdots\subset B_j\subset B_{j+1}\subset\cdots\subset\Omega_T.
\end{equation}
Define $B:=\bigcup_{j=1}^{\infty}B_j$, then $|\Omega_T\backslash B|=0$. Define $P_j:=\{(x,t)\in\Omega_T:|1-u^2|\geq\delta_j\}$. Then
\begin{equation}
\label{3-18}
P_1\subset P_2\subset\cdots \subset P_j\subset P_{j+1}\subset\cdots\subset\Omega_T,
\end{equation}
and $\bigcup_{j=1}^{\infty}P_j=P$. Then, for each $j$, $B_j$ can be split into two parts:
$$
D_j:=B_j\bigcap P_j,\quad\hbox{where}~|1-u^2|>\delta_j~~\hbox{and}~u_i\rightarrow u~\hbox{uniformly},
$$
$$
\widehat{D}_j:=B_j\backslash P_j,\quad\hbox{where}~|1-u^2|\leq\delta_j~~\hbox{and}~u_i\rightarrow u~\hbox{uniformly}.
$$
It then follows from (\ref{3-17}) and (\ref{3-18}) that
\begin{equation}
\label{3-19} D_1\subset D_2\subset\cdots\subset D_j\subset D_{j+1}\subset\cdots\subset D:=B\bigcap P,
\end{equation}
and indeed $D=\bigcup_{j=1}^{\infty}D_j$. Then for any $\psi\in L^2(0,T;L^{\frac{2n}{n-2}}(\Omega_T;\mathbb{R}^n))$, we have
\begin{equation}
\begin{aligned}&
\int_{\Omega_T}M_i(u_i)\nabla\mu_i\cdot\psi dxdt
\\
=&\int_{\Omega_T\backslash B_j}M_i(u_i)\nabla\mu_i\cdot\psi dxdt+\int_{B_j\backslash P_j}M_i(u_i)\nabla\mu_i\cdot\psi dxdt+\int_{B_j\bigcap P_j
}M_i(u_i)\nabla\mu_i\cdot\psi dxdt
\\
\leq&C.
\end{aligned}\label{3-20}\end{equation}
Note that $\{B_j\bigcap P_j\}_{j=1}^{\infty}$ is an increasing sequence of sets with a limit $B\bigcap P$, we can check that $\zeta_{j-1}=\zeta_j$ almost everywhere in $B_{j-1}\bigcap P_{j-1}$. Moreover, we can extend $\zeta_j\in L^2(B_j\bigcap P_j)$ into a function $\tilde{\zeta}\in L^2(B\bigcap P)$ by
\begin{equation} \nonumber
\tilde{\zeta}:=\left\{ \begin{aligned}
        &\zeta_j,\quad\forall x\in B_j\bigcap P_j ,\\
              & 0,\quad\forall x\in (B\bigcap P)\backslash (B_j\bigcap P_j).
                          \end{aligned} \right.
                          \end{equation}
Then, for almost every $x\in B\bigcap P$, there exists a limit of $\tilde{\zeta}_j(x)$ as $j\rightarrow\infty$. Write
$$
\zeta(x)=\lim_{j\rightarrow\infty}\tilde{\zeta}_j(x)\quad\hbox{almost~everywhere~in}~B\bigcap P.
$$
Clearly, $\zeta(x)=\zeta_j(x)$ almost everywhere $x\in B_j\bigcap P_j)$ for all $j$. By a standard diagonal argument, we extract a subsequence such that
\begin{equation}
\label{3-21}
\nabla\mu_{k,N_k}\rightharpoonup \zeta\quad\hbox{weakly~in}~L^2(B_j\bigcap P_j)\quad\hbox{for~all}~j.
\end{equation}
It then follows from (\ref{3-13}) that
\begin{equation}\label{3-22}
\chi_{B_j\bigcap P_j}\sqrt{M_{k,N_k}(u_{k,N_k})}\nabla\mu_{k,N_k}\rightharpoonup \chi_{B_j\bigcap P_j}\sqrt{M(u)}\xi~~\hbox{weakly~in}~L^2(0,T;L^{\frac{2n}{2+n}}(\Omega))~\hbox{for~all}~j.
\end{equation}
Here $\chi_{B_j\bigcap P_j}$ is the characteristic function of $B_j\bigcap P_j\subset\Omega_T$.  By (\ref{3-14}), we know $\sqrt{M_i(u_i))}\rightharpoonup\xi$ weakly in $L^2(\Omega_T)$, we see in fact that $\xi=\sqrt{M(u) }\zeta$ in every $B_j\bigcap P_j$, and hence
\begin{equation}\label{3-23}
\xi=\sqrt{M(u)}\zeta\quad\hbox{in}~B\bigcap P.
\end{equation}
Consequently, by (\ref{3-15}),
$$
\chi_{B\bigcap P}M_{k,N_k}(u_{k,N_k})\nabla\mu_{k,N_k}\rightharpoonup\chi_{B\bigcap P}M(u)\zeta,
$$
weakly in $L^2(0,T;L^{\frac{2n}{n+2}}(\Omega))$.
Replacing $u_i$ by the above mentioned subsequence $u_{k,N_k}$ in (\ref{3-20}), taking limits first as $k\rightarrow\infty$ and then as $j\rightarrow\infty$, we obtain that for every $\psi\in L^2(0,T;L^{\frac{2n}{n+2}}(\Omega;\mathbb{R}^n))$
\begin{equation}\label{3-24}
\int_{\Omega_T}\sqrt{M(u)}\xi\cdot\psi dxdt=\lim_{j\rightarrow\infty}\int_{B_j\bigcap P_j}M(u)\zeta\cdot\psi dxdt=\int_{B\bigcap P}M(u)\zeta\cdot\psi dxdt.
\end{equation}
Compared with (\ref{3-16}), we find that $u$ and $\zeta$ solve the following weak equation
\begin{equation}\label{3-25}\begin{aligned}
\int_0^T\langle\partial_tu,\rho\rangle_{((H^2(\Omega))',H^2(\Omega))}dt=-\int_{B\bigcap P}M(u)\zeta\cdot\nabla\rho dxdt+\beta\cdot\int_0^T\!\!\!\int_{\Omega}\psi'(u)\rho\nabla u dxdt,\end{aligned}
\end{equation}
for all $\rho\in L^2(0,T;H^2(\Omega))$.

The desired relation between $\zeta$ and $u$ is $$\zeta=-\gamma\nabla\Delta u+\varphi'(u)\nabla u.$$ This is a delicate question to be studied here. In fact, while $\varphi'(u)\nabla u\in L^{\infty}(0,T;L^{\frac{2n}{n+2}}(\Omega))$,
the term $\nabla\Delta u$ is only defined in the sense of distributions and may not even be a function, given the known regularity $u\in L^2(0,T;H^1(\Omega))$.

\begin{claim}
\label{claim1}
If for some $j$, the interior of $B_j\bigcap P_j$, denoted by $(B_j\bigcap P_j)^{\circ}$, is not empty, then
$$
\nabla\Delta u\in L^{\frac{2n}{n+2}}((B_j\bigcap P_j)^{\circ}),
$$
and
\begin{equation}\begin{aligned}
\zeta=-\gamma\nabla\Delta u+\varphi'(u)\nabla u,\quad\hbox{in}~(B_j\bigcap P_j)^{\circ}.
\end{aligned}\nonumber\end{equation}
\end{claim}

We also need to consider $\zeta$ is not defined in $\Omega\setminus (B\bigcap P)$.
Note that the value of $\zeta$ in $\Omega_T\setminus (B\bigcap P)$ does not matter as it does not appear in the right hand side of (\ref{3-25}). This ambiguity can be removed in every open subset of $\Omega_T$ in which $\nabla\Delta^2u$ has enough regularity.

\begin{claim}
\label{claim2}
For any open set $U\subset\Omega_T$ in which $\nabla\Delta u\in L^q(U)$ for some $q>1$, where $q$ may depend on $U$, we have
\begin{equation}\begin{aligned}
\zeta=-\gamma\nabla\Delta u+\varphi'(u)\nabla u,\quad\hbox{in}~U.
\end{aligned}\nonumber\end{equation}
\end{claim}

Since the proofs of Claims \ref{claim1}-\ref{claim2} are similar to Subsection 3.2 in \cite{D1}, we omit them here.

Define
$$
\tilde{\Omega}_T:=\bigcup\{U\subset\Omega_T:\nabla\Delta u\in L^p(U)~~\hbox{for~some}~p>1,~~p~\hbox{depending~on}~U\}.
$$
Then $\tilde{\Omega}_T$ is open and
\begin{equation}\begin{aligned}
\zeta=-\gamma\nabla\Delta u+\varphi'(u)\nabla u,\quad\hbox{in}~\tilde{\Omega}_T.
\end{aligned}\nonumber\end{equation}

Now, $\zeta$ is defined in $\tilde{\Omega}_T\bigcup(B\bigcap P)$. In order to extend the definition of $\zeta$ to $\Omega_T$, notice that
$$
\Omega_T\setminus((B\bigcap P)\bigcup\tilde{\Omega}_T)\subset(\Omega_T\setminus P)\bigcup(\Omega_T\setminus B).
$$
Since $|\Omega_T\setminus B|=0$ and $M(u)=0$ in $\Omega_T\setminus P$, the value of $\zeta$ outside of $(B\bigcap P)\bigcup\tilde{\Omega}_T$ does not contribute to the integral on the right hand side of (\ref{3-25}). So, we may just set $\zeta=0$ out side of $(B\bigcap P)\bigcup \tilde{\Omega}_T$.


\subsection*{Acknowledgements}

This paper is supported by the Natural Science Foundation of China for Young Scholar (grant
No. 11401258), Natural Science Foundation of Jiangsu Province for Young Scholar (grant
No. BK20140130) and China Postdoctoral Science Foundation (grant No. 2015M581689
).

\end{document}